\documentclass[journal]{IEEEtran}

\usepackage{amsmath}
\usepackage{amsfonts}
\usepackage{tensor}
\usepackage{microtype}
\usepackage{nomencl}
\usepackage{graphicx}
\usepackage{subfigure}
\usepackage{booktabs}
\usepackage{multirow}
\usepackage{mathtools,cases}
\usepackage{xcolor}
\usepackage{hyperref}
\usepackage{etoolbox} 

\makenomenclature

\renewcommand\nomgroup[1]{\vspace{3mm}
  \item[\bfseries
  \ifstrequal{#1}{S}{Sets}{
  \ifstrequal{#1}{V}{\hspace*{-2mm}Variables}{
  \ifstrequal{#1}{F}{\hspace*{-4mm} Functions and functionals}{}}}]
\hspace*{-\leftmargin}\vspace{3mm}
}

\begin{document}

\urlstyle{tt}

\title{A General Stochastic Optimization Framework for Convergence Bidding}
%
%
%

\author{Letif~Mones
        and~Sean~Lovett
\thanks{All authors are with Invenia Labs, 95 Regent Street, Cambridge, CB2 1AW, United Kingdom (e-mails: \{firstname.lastname\}@invenialabs.co.uk).}
}

\maketitle

\begin{abstract}
Convergence (virtual) bidding is an important part of two-settlement electric power markets as it can effectively reduce discrepancies between the day-ahead and real-time markets.
Consequently, there is extensive research into the bidding strategies of virtual participants aiming to obtain optimal bids to submit to the day-ahead market.
In this paper, we introduce a price-based general stochastic optimization framework to obtain optimal convergence bid curves.
Within this framework, we develop a computationally tractable linear programming-based optimization model, which produces bid prices and volumes simultaneously.
We also show that different approximations and simplifications in the general model lead naturally to state-of-the-art convergence bidding approaches, such as self-scheduling and opportunistic approaches.
Our general framework also provides a straightforward way to compare the performance of these models, which is demonstrated by numerical experiments on the California (CAISO) market.
\end{abstract}

\begin{IEEEkeywords}
Convergence bidding, virtual bidding, electricity market, stochastic optimization, portfolio optimization.
\end{IEEEkeywords}

\IEEEpeerreviewmaketitle


\nomenclature[S]{$\mathcal{P}$}{Convex set of supply and demand prices.}
\nomenclature[S]{$\mathcal{W}$}{Convex set of supply and demand volumes.}
\nomenclature[S]{$\mathcal{B}$}{Convex set of bids (supply and demand prices and volumes).}

\nomenclature[V]{$F$}{Number of physical participants in the day-ahead market.}
\nomenclature[V]{$N$}{Number of biddable nodes in the target period.}
\nomenclature[V]{$T$}{Number of training samples.}
\nomenclature[V]{$S$}{Number of maximum allowed segments in bid curve.}
\nomenclature[V]{$V$}{Number of all other virtual participants in the day-ahead market.}
\nomenclature[V]{$W$}{Total volume distributed among nodes and positions in stochastic optimization problem.}
\nomenclature[V]{$\overline{W}, \underline{W}$}{Upper and lower bounds on net volume in stochastic optimization problem.}
\nomenclature[V]{$X, x^{*}$}{Training feature matrix and test feature vector for a probabilistic model.}
\nomenclature[V]{$Y$}{Training price matrix for a probabilistic model.}
\nomenclature[V]{$\overline{P}, \underline{P}$}{Number of potential bid prices for supply and demand positions.}
\nomenclature[V]{$K$}{Number of samples used to compute expected shortfall.}
\nomenclature[V]{$\lambda^{n}_{t}, \pi^{n}_{t}$}{Day-ahead and real-time LMPs at node $n$ at index $t$ in the training data.}
\nomenclature[V]{$\Delta^{n}_{t}$}{Delta (spread) LMP at node $n$ at index $t$ in the training data.}
\nomenclature[V]{$r, \tilde{r}$}{Hourly revenue and volume-normalized revenue variables.}
\nomenclature[V]{$\rho, \tilde{\rho}$}{Risk and volume-normalized risk upper bound parameters.}
\nomenclature[V]{$\tau, z$}{Variables for computing sample-based expected shortfall.}
\nomenclature[V]{$\theta$}{Parameter vector of a probabilistic model.}
\nomenclature[V]{$\mu, \Sigma$}{Mean vector and covariance matrix of predicted delta price distribution.}
\nomenclature[V]{$b$}{Submitted bid including $(\overline{w}, \underline{w}, \overline{p}, \underline{p})$.}
\nomenclature[V]{$\overline{p}, \underline{p}$}{Supply and demand price variables.}
\nomenclature[V]{$\overline{w}, \underline{w}$}{Supply and demand volume variables.}
\nomenclature[V]{$\overline{v}, \underline{v}$}{Cleared supply and demand volume variables.}
\nomenclature[V]{$\overline{d}, \underline{d}$}{Minimum distance in prices between supply and demand bid segments.}
\nomenclature[V]{$z_{\lambda}$, $z_{\pi}$}{Parameters of the day-ahead and real-time OPF problems.}

\nomenclature[F]{$\varphi, \psi$}{Abstract utility functions used in stochastic optimizations.}
\nomenclature[F]{$\mathbb{F}$, $\mathbb{G}$}{Abstract functionals used in stochastic optimizations.}
\nomenclature[F]{$\mathbb{E}$}{Expected value.}
\nomenclature[F]{$\mathbb{V}$}{Variance.}
\nomenclature[F]{$\mathbb{ES}^{\alpha}$}{Expected shortfall (conditional value at risk) with quantile parameter $\alpha$.}
\nomenclature[F]{$\mathcal{L}, \mathcal{D}$}{Abstract deterministic functions to compute day-ahead and delta prices.}

\printnomenclature

\section{Introduction}
\label{sec:intro}

\IEEEPARstart{C}{onvergence} or virtual bidding plays a critical role in most electric power markets in the United States.
Independent System Operators (ISOs) use a two-settlement approach to aid price convergence (i.e. reduce the price gap) between the day-ahead and real-time (spot) markets.
Under this mechanism, participants can submit convergence supply (increment) and demand (decrement) bids at each target hour to the day-ahead market.
Convergence bids are input into the market clearing calculations in the same way as physical bids, but unlike physical bids, they are settled financially and carry no obligation to supply or consume physical power in real time.
The revenue of the participant depends on the difference between day-ahead and real-time Locational Marginal Prices (LMPs) \cite{Liu09}.
There is a wide consensus that virtual bidding improves the efficiency of electricity market operation by converging prices~\cite{Isemonger06, Celebi10, Jha13, Li15, Woo15, Tang16, Hogan16, Mather17} and increasing day-ahead unit commitment efficiency~\cite{Kazempour18a, Kazempour18b}, hedges financial risks~\cite{Celebi10, Hogan16}, increases market liquidity~\cite{Hogan16}, reduces price volatility~\cite{Hadsell07, Woo15}, and helps prevent the exercise of market power~\cite{Isemonger06, Hogan16}.

The problem of how to select optimal convergence bids is important both to virtual participants and to the efficient operation of electricity markets.
While this problem bears some relationship to classical problems of risk-averse portfolio optimization, it has several unique aspects.
The most important is that the clearance of a convergence bid is not based on a simple bilateral auction, but rather on the solution of the day-ahead unit commitment and optimal power flow (OPF) problems.
A convergence bid on a given node at a given target hour is submitted as a volume-price pair (i.e. bid volume of power in unit of MWh and bid price in unit of \$/MWh).
In the day-ahead market, a convergence bid with a supply position is cleared if the (day-ahead) LMP at the specific hour and node is greater or equal than the bid price and similarly, a convergence bid with demand position is cleared if the corresponding LMP is less or equal than the bid price.
As part of the optimization problem, the cleared convergence bid will alter the day-ahead market solution and therefore the day-ahead LMPs.
Accordingly, a convergence bidding approach can be modelled by a bilevel optimization problem, where the upper level is a financial problem subject to lower level market clearance problems~\cite{Mehdipourpicha20, Xiao21, Kohansal21}.

Ideally, virtual participants want to obtain nodal bid curves including segments with optimal bid volumes and prices.
This requires solving stochastic optimization problems that treat both bid volumes and prices as decision variables at the same time.
However, even without the bilevel consideration, co-optimizing both prices and volumes results in large mixed-integer optimization problems, which are typically computationally infeasible for realistic data.
This leads to the consideration of simplified, computationally feasible versions of this problem.
A typical approach is to focus only on a subset of the original variables in the decision making process, while using fixed values or simple heuristics to obtain the rest.
For example, optimal bid volumes can be derived with fixed prices~\cite{Wang19, Li22}, or optimal prices with fixed quantity of power~\cite{Baltaoglu19, Saez-Gallego18, Samani21}.

In this paper, we consider a general framework for financially-based convergence bidding and derive an optimization model that provides optimal virtual bid curves.
More specifically, our contributions are as follows:
\begin{itemize}
    \item We introduce a fully price-based stochastic optimization framework for obtaining optimal bid curves.
    \item The general framework is primarily based on some predictive joint distribution of the day-ahead and delta prices at the target hours and includes both bid volumes and bid prices as decision variables (referred to as a VP-model).
    \item We provide a detailed discussion on approaches for approximating the joint distribution of the day-ahead and delta prices.
    \item We show that the general framework encompasses many state-of-the-art approaches as special cases. More precisely, omitting bid prices from the decision making model and focusing only on bid volumes leads to \emph{self-scheduling} approaches or V-models, while omitting bid volumes and focusing only on bid prices leads to \emph{opportunistic} approaches or P-models.
    \item Within the general framework, we introduce a new optimization model that 1. provides optimal bid curves (i.e. with multiple bid segments) 2. for both supply and demand positions 3. by finding optimal bid prices and volumes simultaneously 4. using a linear programming problem. The model avoids integer variables and is therefore computationally affordable for large numbers of nodes and training data.
    \item The general framework also provides a straightforward way to carry out a fair comparison among different (VP, V and P) models. We demonstrate it by comparing state-of-the-art (V and P) approaches to our model (VP) on the California (CAISO) market.
\end{itemize}

\section{Methodology}

\subsection{Setup}

Suppose we have an electricity market with $N$ nodes potentially biddable by virtual participants. 
Our aim is to derive bids (volumes and prices) for the target period of the day-ahead market using historical price data.
Let $\lambda_{t}^{n}$, $\pi_{t}^{n}$ and $\Delta_{t}^{n} = \lambda_{t}^{n} - \pi_{t}^{n}$ be the corresponding locational marginal prices (LMPs) of the day-ahead and real-time markets as well as the delta (spread) prices, where $n=1, \dots, N$ and $t=1, \dots, T$, with $T$ number of training samples.
Let $\lambda_{t}$, $\pi_{t}$ and $\Delta_{t}$ denote the vectors of all corresponding nodal prices at time $t$.
Furthermore, without any indices, $\Delta$ and $\lambda$ represent vector random variables with dimension $N$, while $\Delta^{n}$ and $\lambda^{n}$ denote their $n$th components, respectively.  

Electric power markets support the submission of convergence bid curves, i.e. multiple bid segments with different bid prices for the same target hour and node.
Therefore, we wish to derive optimal supply and demand convergence bid curves, each with a maximum number of $S$ segments including supply volumes ${\overline{w}}^{n}_{s} \ge 0$ and prices ${\overline{p}}^{n}_{s}$ as well as demand volumes ${\underline{w}}^{n}_{s} \le 0$ and prices ${\underline{p}}^{n}_{s}$, where $n=1, \dots, N$ and $s=1, \dots, S$.
In this work, we follow a `block' interpretation of bid curves, where a bid curve consists simply of several bids, which are cleared independently (cf. MISO~\cite{MISO-BPM}).
This can be converted straightforwardly to the `tiered` interpretation used by some markets, where a bid curve consists of monotonically increasing segments, at most one of which is cleared (cf. CAISO~\cite{CAISO-BPM}).

We emphasize that in the general framework, at the target hour for a given node $n$ we can have the following four possible bid types, where $s$ and $s'$ denote bid segment indexes:
\begin{itemize}
    \item No bid position: $\forall s: {\overline{w}}^{n}_{s} = {\underline{w}}^{n}_{s} = 0$.
    \item Supply position only: $\exists s: {\overline{w}}^{n}_{s} > 0$ and $\forall s': {\underline{w}}^{n}_{s'} = 0$.
    \item Demand position only: $\forall s: {\overline{w}}^{n}_{s} = 0$ and $\exists s': {\underline{w}}^{n}_{s'} < 0$.
    \item Both supply and demand positions: $\exists s, s': {\overline{w}}^{n}_{s} > 0$ and ${\underline{w}}^{n}_{s'} < 0$.
\end{itemize}

\subsection{General stochastic optimization problem for convergence bidding}
\label{sec:general_so}

We suggest the following fully price-based general form of stochastic optimization problem for convergence bidding, in order for any specific virtual participant to determine optimal supply and demand bid curves with corresponding volumes and prices:
\begin{equation}
    \begin{aligned}
        \max_{\substack{{\overline{w}, \overline{p} \in \mathbb{R}^{S \times N}} \\ {\underline{w}, \underline{p} \in \mathbb{R}^{S \times N}}}} & \ \mathbb{F}_{p(\Delta, \lambda)} \left[\varphi(\Delta, \lambda, \overline{w}, \underline{w}, \overline{p}, \underline{p})\right] \\
        \mathrm{s.t.} & \ \mathbb{G}_{p(\Delta, \lambda)} \left[\psi(\Delta, \lambda, \overline{w}, \underline{w}, \overline{p}, \underline{p})\right] \le 0, \\
        & \ (\overline{w}, \underline{w}) \in \mathcal{W}, \\
        & \ (\overline{p}, \underline{p}) \in \mathcal{P}, \\
    \end{aligned}
    \label{eq:general_so}
\end{equation}
where $p(\Delta, \lambda)$ is a given joint distribution of the random variables representing delta and day-ahead prices.
$\varphi$ and $\psi$ are some utility functions of these prices (e.g. financial revenue), as well as volumes ($\overline{w}$, $\underline{w}$) and bid prices ($\overline{p}$, $\underline{p}$), which are the decision variables.
Clearly, $\varphi$ and $\psi$ are also random variables due to $\Delta$ and $\lambda$.
In order to take uncertainties of the outcomes of the decision into account, statistical functionals $\mathbb{F}$ (e.g. expectation value) and $\mathbb{G}$ (e.g. some risk measure) of the utility functions are used.
$\mathcal{W}$ and $\mathcal{P}$ are (convex) sets of the volumes and bid prices, respectively.

We note that in a more general framing, the delta and day-ahead prices are deterministic functions of the bids of all physical and virtual participants based on the corresponding optimal power flow calculations carried out by the ISOs:
\begin{equation}
    \begin{aligned}
        \max_{\substack{b}} & \ \mathbb{F}_{p(b, \{ b^{\mathrm{phys}}_{i} \}_{i=1}^{F}, \{ b^{\mathrm{virt}}_{j} \}_{j=1}^{V}, z_{\lambda}, z_{\pi})} \left[\varphi(\Delta, \lambda, b)\right] \\
        \mathrm{s.t.} & \ \mathbb{G}_{p(b, \{ b^{\mathrm{phys}}_{i} \}_{i=1}^{F}, \{ b^{\mathrm{virt}}_{j} \}_{j=1}^{V}, z_{\lambda}, z_{\pi})} \left[\psi(\Delta, \lambda, b)\right] \le 0, \\
        & \ b \in \mathcal{B}, \\
        \mathrm{with\ } & \lambda = \mathcal{L} \left(b, \{ b^{\mathrm{phys}}_{i} \}_{i=1}^{F}, \{ b^{\mathrm{virt}}_{j} \}_{j=1}^{V}, z_{\lambda} \right), \\
        & \Delta = \mathcal{D} \left(\lambda, z_{\pi} \right), \\
    \end{aligned}
    \label{eq:very_general_so}
\end{equation}
where $b = (\overline{w}, \underline{w}, \overline{p}, \underline{p})$ is the full set of bids submitted by our specific virtual participant and $\mathcal{B}$ is a convex set, while $b^{\mathrm{phys}}_{i} = (\overline{w}^{\mathrm{phys}}_{i}, \underline{w}^{\mathrm{phys}}_{i}, \overline{p}^{\mathrm{phys}}_{i}, \underline{p}^{\mathrm{phys}}_{i})$, $i=1, \dots, F$ and $b^{\mathrm{virt}}_{j} = (\overline{w}^{\mathrm{virt}}_{j}, \underline{w}^{\mathrm{virt}}_{j}, \overline{p}^{\mathrm{virt}}_{j}, \underline{p}^{\mathrm{virt}}_{j})$, $j=1, \dots, V$ represent bids from physical and all other virtual participants and $F$ and $V$ are the number of physical and all other virtual participants, respectively.
$\mathcal{L}$ is a deterministic function of all bids and some parameters $z_{\lambda}$ to provide day-ahead LMP prices, while $\mathcal{D}$ is a deterministic function of the day-ahead prices and some parameters $z_{\pi}$ associated to the real-time market.
In this case, the uncertainty in the corresponding stochastic problem arises due to the distribution $p(\{b^{\mathrm{phys}}_{i} \}_{i=1}^{F}, \{b^{\mathrm{virt}}_{j} \}_{j=1}^{V}, z_{\lambda}, z_{\pi})$.
The potential advantage of this approach is that the price impact of bids submitted by our specific virtual participant is explicitly considered.
Also, in principle, this approach could provide accurate estimates of the dependency between LMPs (both $\lambda$ and $\Delta$) by incorporating structure from the corresponding OPF solutions.
However, there are several challenges with this most general model that make it impractical to obtain optimal virtual bids for a specific virtual participant.
First, it would require models for computing day-ahead and real-time prices that are adequate approximations of the corresponding OPF and LMP calculations of the ISO; however in practice, full details and parameters ($z_{\lambda}$ and $z_{\pi}$) of the OPF formulations are typically unavailable to the participant at time of bidding.
Since these models themselves are optimization problems, the resulting stochastic optimization problem is rather a complicated one (i.e. non-linear or even non-convex, mixed-integer and possibly nested optimization problems).
Second, it would also need an appropriate distribution of bids by the other participants and an extensive sampling of this distribution that would significantly increase the size and/or the number of the optimization problems to solve.
One way to approximate the most general model is to apply bilevel formulations (\cite{Mehdipourpicha20, Xiao21, Kohansal21}) that try to address the clearing mechanism by solving some simplified OPF on the lower level.
These problems, among other extensions, could be considered by using $p(\Delta, \lambda\ |\ \overline{w}, \underline{w}, \overline{p}, \underline{p})$, that is, accounting only for the impact of actions of the specific virtual participant on the joint distribution of prices.
Problem~\ref{eq:general_so} is a further simplification, where the clearing process is not modelled explicitly and uncertainty about the behaviour of all (physical and virtual) participants is incorporated into $p(\Delta, \lambda)$. 

A typical choice for the functions $\varphi$ and $\psi$ in problem~\ref{eq:general_so} is the financial revenue.
If cleared, a supply bid segment at node $n$ results in a payment of $\overline{w}^{n}_{s} \Delta^{n}_{*}$ if $\lambda^{n}_{*} \ge \overline{p}^{n}_{s}$ and similarly, a demand bid segment has a payment of $\underline{w}^{n}_{s} \Delta^{n}_{*}$ if $\lambda^{n}_{*} \le \underline{p}^{n}_{s}$, where the subscript $*$ denotes the corresponding target hour.
However, our aim is to take situations into account when the bid prices are set such that certain bids or bid segments do not clear (i.e. their contribution is 0).
Therefore, we consider the following expression: 
\begin{equation}
    \begin{aligned}
        \varphi(\Delta, \lambda, \overline{w}, \underline{w}, \overline{p}, \underline{p}) & = \psi(\Delta, \lambda, \overline{w}, \underline{w}, \overline{p}, \underline{p}) \\
        & = \left(\sum \limits_{s=1}^{S} 1_{\lambda \succeq {\overline{p}}_{s}} \odot {\overline{w}}_{s} + 1_{\lambda \preceq {\underline{p}}_{s}} \odot {\underline{w}}_{s} \right)^{T} \Delta,
    \end{aligned}
\end{equation}
where $\odot$ designates the Hadamard (i.e. pointwise) product between two vectors.
$1_{\lambda \succeq {\overline{p}}_{s}}$ and $1_{\lambda \preceq {\underline{p}}_{s}}$ denote vector-valued indicator functions, whose $n$th component is 1 if the corresponding position is cleared, i.e. a supply bid segment is cleared if its bid price is less than or equal to the day-ahead price ($\lambda^{n} \ge {\overline{p}}^{n}_{s}$) and similarly, a demand bid segment is cleared if its bid price is greater than or equal to the day-ahead price ($\lambda^{n} \le {\underline{p}}^{n}_{s}$).

Using the above definition of the revenue results in the following stochastic optimization problem:
\begin{equation}
    \begin{aligned}
        \max_{\substack{{\overline{w}, \overline{p} \in \mathbb{R}^{S \times N}} \\ {\underline{w}, \underline{p} \in \mathbb{R}^{S \times N}}}} & \ \mathbb{F}_{p(\Delta, \lambda)}\left[\left(\sum \limits_{s=1}^{S} 1_{\lambda \succeq {\overline{p}}_{s}} \odot {\overline{w}}_{s} + 1_{\lambda \preceq {\underline{p}}_{s}} \odot {\underline{w}}_{s} \right)^{T} \Delta\right] \\
        \mathrm{s.t.} & \ \mathbb{G}_{p(\Delta, \lambda)}\left[\left(\sum \limits_{s=1}^{S} 1_{\lambda \succeq {\overline{p}}_{s}} \odot {\overline{w}}_{s} + 1_{\lambda \preceq {\underline{p}}_{s}} \odot {\underline{w}}_{s} \right)^{T} \Delta \right] \le 0, \\
        & \ (\overline{w}, \underline{w}) \in \mathcal{W}, \\
        & \ (\overline{p}, \underline{p}) \in \mathcal{P}. \\
    \end{aligned}
    \label{eq:return_so}
\end{equation}

It is common to choose $\mathbb{F}$ as the expected value of the revenue according to $p(\Delta, \lambda)$.
Also, $\mathbb{G}$ is typically some risk measure (or function of risk measure) of the revenue according to $p(\Delta, \lambda)$, such as variance, value at risk \cite{Jorion00} or expected shortfall (conditional value at risk) \cite{Rockafellar00}.
We also note that problem~\ref{eq:return_so} represents other formulations as well.
For instance, one could use a trade-off objective (risk-adjusted revenue) for $\mathbb{F}$ with an empty $\mathbb{G}$.
Similarly, setting $\mathbb{F}$ to the negative risk and $\mathbb{G}$ to some function of the negative expected revenue leads to a risk minimization with a lower bound on the expected revenue.

As a concrete example, we show the expected revenue maximization with an upper bound on the expected shortfall using single segment bids (i.e. $S=1$):
\begin{equation}
    \begin{aligned}
        \max_{\substack{{\overline{w}, \overline{p} \in \mathbb{R}^{N}} \\ {\underline{w}, \underline{p} \in \mathbb{R}^{N}}}} & \ \mathbb{E}_{p(\Delta, \lambda)}\left[(1_{\lambda \succeq \overline{p}} \odot \overline{w} + 1_{\lambda \preceq \underline{p}} \odot   \underline{w})^{T} \Delta \right] \\
        \mathrm{s.t.} & \ \mathbb{ES}^{\alpha}_{p(\Delta, \lambda)}\left[(1_{\lambda \succeq \overline{p}} \odot \overline{w} + 1_{\lambda \preceq \underline{p}} \odot   \underline{w})^{T} \Delta \right] \le \rho, \\
        & \ (\overline{w}, \underline{w}) \in \mathcal{W}, \\
        & \ (\overline{p}, \underline{p}) \in \mathcal{P}, \\
    \end{aligned}
    \label{eq:return_es_so}
\end{equation}
where $\alpha$ is the quantile parameter of the expected shortfall and $\rho$ is the upper bound.

A set of typical constraints that might define $\mathcal{W}$ is:
\begin{itemize}
    \item Sign of the supply and demand volume components: $\overline{w} \succeq 0$ and $\underline{w} \preceq 0$. By convention, this set of constraints will always be present.    
    \item Lower and upper bounds on volume components: $\overline{w}_{\mathrm{min}} \preceq \overline{w} \preceq \overline{w}_{\mathrm{max}}$ and $\underline{w}_{\mathrm{min}} \preceq \underline{w} \preceq \underline{w}_{\mathrm{max}}$.
    \item Total attempted volume (L1-norm): $\| \overline{w} - \underline{w}\|_{1, 1} \le W$, where $W > 0$ is the total volume and $\| \cdot \|_{1, 1}$ is the entry-wise matrix norm.
    \item Net volume constraint: $\underline{W} \le \| \overline{w}\|_{1, 1} - \| \underline{w} \|_{1, 1} \le \overline{W}$, where $\underline{W}$ and $\overline{W}$ are the lower and upper bounds on the net volume. A specific case of this constraint is net-zero convergence bidding (i.e. $\underline{W} = \overline{W} = 0$ for each hour)~\cite{Kohansal21}.
    \item Additional specific constraints for transaction fees, uplift costs and sensitivity of convergence bids~\cite{Li22}, which could also be approximated by piecewise-linear functions.
\end{itemize}

For defining $\mathcal{P}$, the main relevant constraints might be:
\begin{itemize}
    \item Lower and upper bounds on price components: $\overline{p}_{\mathrm{min}} \preceq \overline{p} \preceq \overline{p}_{\mathrm{max}}$ and $\underline{p}_{\mathrm{min}} \preceq \underline{p} \preceq \underline{p}_{\mathrm{max}}$.
    \item Minimum distance between adjacent bid segment prices: $|\overline{p}^{n}_{s} - \overline{p}^{n}_{s'}| \ge \overline{d}^{n}$ and $|\underline{p}^{n}_{s} - \underline{p}^{n}_{s'}| \ge \underline{d}^{n}$ for $i=1,\dots,N$, $s, s' = 1,\dots, S$ and $s \ne s'$.
\end{itemize}

Finally, we note that the target period in problem~\ref{eq:general_so} could be hourly or daily.
In the case of an hourly approach, for each target hour, a separate and independent optimization problem would be defined and solved, while for a daily approach, bids would be derived from a single optimization problem for all target hours of the day (conceptually, the same node at different target hours of the day is treated as distinct nodes).

\subsection{Joint delta and day-ahead distribution}
\label{sec:joint}
The key quantity in the above general model is the joint distribution of the delta and day-ahead prices for the target hour: $p(\Delta, \lambda)$.
There are three main approaches that could be used to approximate this distribution:

\begin{itemize}
    \item Using a probabilistic model: $p(\Delta, \lambda | x^{*}, X, Y, \theta)$, where $x^{*}$ and $X = \left\{ x_{t} \right\}_{t=1}^{T}$ denote the test (i.e. at target hour) and training features, $Y = \left\{ \Delta_{t}, \lambda_{t} \right\}_{t=1}^{T}$ are the training prices and $\theta$ are model parameters. The parameters $\theta$ could be optimized for example by maximum likelihood, or by Bayesian inference. The advantage of this approach is that given an accurate inductive bias of the model, limited training data is needed to obtain the joint distribution, minimizing errors caused by non-stationarity, i.e. the fact that $p(\Delta, \lambda)$ changes in time. However, one limitation of this approach is that in practice usually a simple multivariate distribution (e.g. some elliptical distribution) is used in order to obtain an analytically closed-form stochastic optimization model. Such simple distributions will struggle to model features such as multimodality and heavy-tailedness of the underlying true distribution. Another challenge is to find and validate an appropriate inductive bias.
    \item Using the training samples $\left\{ \Delta_{t}, \lambda_{t} \right\}_{t=1}^{T}$ directly, which naturally models any complicated distribution, although with finite samples. However, in order to compute sufficient statistics, this approach requires a rather large amount of training data, and is therefore more exposed to the effects of non-stationarity. Since the number of variables in the optimization problem is proportional to the training data size, it is advisable to use a linear programming based model for this approach. We note that a marginalised version of this approach was used in~\cite{Baltaoglu19, Saez-Gallego18, Samani21}, representing individual joint distributions per node (i.e. $\left\{ \Delta_{t}^{n}, \lambda_{t}^{n} \right\}_{t=1}^{T}$).
    \item Finally, we note that a combination of the above two approaches is also possible when one uses a probabilistic model with a more sophisticated distribution. Using the distribution directly within the stochastic optimization problem would not necessarily lead to an analytically closed-form expression; however, the distribution can be sampled and these samples can be used. This approach leads to similar optimization problems as using empirical samples directly.
\end{itemize}

\begin{figure}[!ht]
    \centerline{\includegraphics[width=\columnwidth]{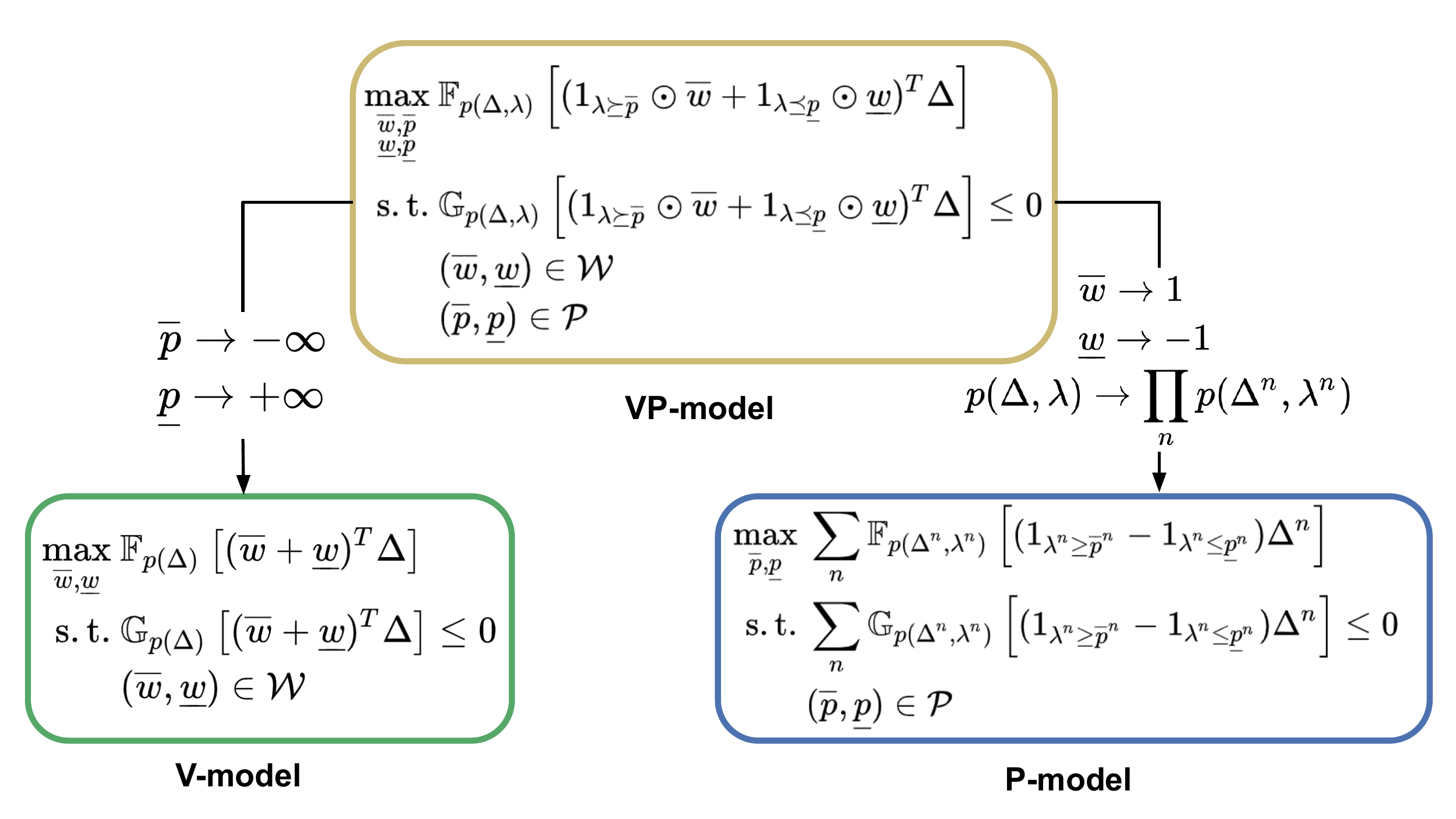}}
    \caption{The general stochastic optimization framework (VP) and its primary simplifications (V, P) using single segment bids.}
    \label{fig:vp_v_p_models}
\end{figure}

\subsection{Volume-price (VP) models}

It is an interesting question whether a parametric distribution for $p(\Delta, \lambda)$ exists for which problem \ref{eq:return_so} has a closed-form expression with appropriate risk measure.
The challenge is that the optimization variables $\overline{p}$ and $\underline{p}$ appear in the indicator function along with $\lambda$. A more straightforward approach, however, is to use sample-based optimization.
A concrete model using the expected shortfall as risk measure (problem~\ref{eq:return_es_so}) can be written as:
\begin{equation}
    \begin{aligned}
        \max_{\substack{{\overline{w}, \overline{p} \in \mathbb{R}^{S \times N}} \\ {\underline{w}, \underline{p} \in \mathbb{R}^{S \times N}} \\ {z \in \mathbb{R}^{T}, \tau \in \mathbb{R}} \\ {\overline{v}, \underline{v} \in \mathbb{R}^{T \times S \times N}}}} & \ \frac{1}{T} \sum \limits_{t=1}^{T} r_{t} \\
        \mathrm{s.t.} & \ -\tau + \frac{1}{K} \sum \limits_{t=1}^{T} z_{t} \le \rho, \\
        & \ z_{t} \ge 0, && \forall t \\
        & \ z_{t} \ge \tau - r_{t}, && \forall t \\
        & \ r_{t} = \left( \sum \limits_{s=1}^{S}\overline{v}_{t, s} + \underline{v}_{t, s} \right)^{T} \Delta_{t}, && \forall t \\
        & \ \overline{v}_{t,s}^{n} = \begin{cases} \overline{w}^{n}_{s} \mathrm{\ if \ } \lambda_{t}^{n} \ge \overline{p}^{n}_{s} \\ 0 \mathrm{\ \ \ else} \\ \end{cases}, && \forall t, \ \forall s, \ \forall n \\
        & \ \underline{v}_{t,s}^{n} = \begin{cases} \underline{w}^{n}_{s} \mathrm{\ if \ } \lambda_{t}^{n} \le \underline{p}^{n}_{s} \\ 0 \mathrm{\ \ \ else} \\ \end{cases}, && \forall t, \ \forall s, \ \forall n \\
        & \ (\overline{w}, \underline{w}) \in \mathcal{W}, \\
        & \ (\overline{p}, \underline{p}) \in \mathcal{P}, \\
    \end{aligned}
    \label{eq:vp_sample}
\end{equation}
where the objective to maximize is the expected revenue computed from the cleared segments.
The first three constraints with variables $z \in \mathbb{R}^{T}$ and $\tau \in \mathbb{R}$ are the corresponding expressions for the sample expected shortfall and its upper bound with $K = \lfloor \alpha T \rfloor$\cite{Rockafellar00}.
Variables $\overline{v}_{t, s}^{n}$ and $\underline{v}_{t, s}^{n}$ express the volumes of the cleared supply and demand of segment $s$ on node $n$ at time $t$.
The above problem leads to a mixed-integer linear programming (MILP) problem, which requires $\mathcal{O}(SN)$ memory for the continuous variables, and a large amount of binary variables to compute $\overline{v}_{t, s}^{n}$ and $\underline{v}_{t, s}^{n}$ (with $\mathcal{O}(TSN)$ memory requirement).
Given that MILP is NP-hard, problem~\ref{eq:vp_sample} with sufficient training data is rather expensive to solve even for small~$N$.

However, the above model can be converted to a pure linear programming (LP) problem by introducing a larger number of (non-binary) variables with a memory requirement of $\mathcal{O}(TN)$ that results in a problem that can be solved within reasonable time.
We emphasize that the conversion is exact and no relaxation is required.
Below, we describe this idea, which we call the sample-based mean-expected shortfall bid curve volume-price model (`sample-VP').

It is clear that in problem \ref{eq:vp_sample} the optimal bid prices for a given node $n$ (i.e. $\overline{p}^{n}_{s}$ and $\underline{p}^{n}_{s}$ for $s = 1, \dots, S$) can be selected without loss of generality from the corresponding training day-ahead prices ($\left\{ \lambda^{n}_{t} \right\}_{t=1}^{T}$ for $n = 1, \dots, N$).
The challenge of the optimization task is to compute the corresponding revenue values correctly for a given bid price, i.e. only the cleared delta prices are considered while the rejected ones are not.
In problem~\ref{eq:vp_sample}, similarly to previous price-based models~\cite{Saez-Gallego18, Samani21}, this is carried out by introducing binary variables that set the volumes of rejected bids to 0 depending on the bid price.

However, the possible cleared and rejected delta prices can be precomputed for all possible bid prices taken from the training day-ahead prices.
Let us consider the $n$th single node for which we wish to find an optimal supply bid price. Depending on $\mathcal{P}$ there are $\overline{P} \le T$ number of possible bid prices from the available day-ahead prices $\overline{p}^{n}_{q} \in \left\{ \lambda_{t}^{n} \right\}_{t=1}^{T}$, $q=1, \dots, \overline{P}$.
We can create a $T \times \overline{P}$ indicator matrix whose columns correspond to the potential bid prices and rows represent the actual indices of the training data: $\overline{I}^{n} = \lambda^{n} \otimes_{\ge} \overline{p}^{n}$, where the symbol $\otimes_{\ge}$ is introduced with the analogy to the outer product of two vectors ($a \otimes b = ab^{T}$), i.e. $a \otimes_{\ge} b = a \succeq b^{T}$.
In other words, $(\overline{I}^{n})_{tq} = 1_{\lambda^{n}_{t} \ge \overline{p}^{n}_{q}}$ shows whether the $t$-th training data point is cleared when the bid price is $\overline{p}^{n}_{q}$.
Using this indicator matrix, the corresponding cleared-rejected supply delta price matrix can be easily constructed: 
\begin{equation}
    \overline{\Delta}^{n} = (\Delta^{n} \otimes 1_{\overline{P}}) \odot \overline{I}^{n} = (\Delta^{n} \otimes 1_{\overline{P}}) \odot (\lambda^{n} \otimes_{\ge} \overline{p}^{n}),
\end{equation}
where $\Delta^{n} \otimes 1_{\overline{P}}$ creates a $T \times \overline{P}$ matrix whose each row is $\Delta^{n}$.
Therefore, the value of $(\overline{\Delta}^{n})_{tq}$ is $\Delta^{n}_{t}$ if cleared and 0 if rejected with $\overline{p}^{n}_{q}$ bid price.

Our next step is to introduce the supply volume vector $0 \preceq \overline{w}^{n} \in \mathbb{R}^{\overline{P}}$ that includes potential volumes for each possible bid price.
The vector $\overline{\Delta}^{n} \overline{w}^{n}$ therefore represents the revenues for each training data point for node $n$.
Note that this approach naturally encompasses bid curves, as multiple non-zero elements are allowed in $\overline{w}^{n}$.
A similar expression can be derived for demand positions as well ($\underline{\Delta}^{n} \underline{w}^{n}$). 
Putting the contributions from all nodes and positions together we have the following optimization problem:
\begin{equation}
    \begin{aligned}
        \max_{\substack{{\overline{w} \in \mathbb{R}^{\overline{P} \times N}} \\ {\underline{w} \in \mathbb{R}^{\underline{P} \times N}} \\ {z \in \mathbb{R}^{T}, \tau \in \mathbb{R}}}} & \ \frac{1}{T} \sum \limits_{t=1}^{T} r_{t} \\
        \mathrm{s.t.} & \ -\tau + \frac{1}{K} \sum \limits_{t=1}^{T} z_{t} \le \rho, \\
        & \ z_{t} \ge 0, & \forall t \\
        & \ z_{t} \ge \tau - r_{t}, & \forall t \\
        & \ r = \sum \limits_{n=1}^{N} \left( \overline{\Delta}^{n} \overline{w}^{n} + \underline{\Delta}^{n} \underline{w}^{n} \right), \\
        & \ (\overline{w}, \underline{w}) \in \mathcal{W}, \\
        \mathrm{with} & \ \overline{\Delta}^{n} = (\Delta^{n} \otimes 1_{\overline{P}}) \odot (\lambda^{n} \otimes_{\ge} \overline{p}^{n}), & \forall n \\
        & \ \underline{\Delta}^{n} = (\Delta^{n} \otimes 1_{\underline{P}}) \odot (\lambda^{n} \otimes_{\le} \underline{p}^{n}), & \forall n \\
        & \ (\overline{p}, \underline{p}) \in \mathcal{P}. \\
    \end{aligned}
    \label{eq:fvp_sample}
\end{equation}
Note that in the above problem, $\overline{\Delta}^{n}, \underline{\Delta}^{n}, \overline{p}^{n}$ and $\underline{p}^{n}$ are precomputed for $n=1,\dots, N$.
Once the above LP is solved for the volumes, the corresponding bid price values can be easily identified from the non-zero weights.

Unlike the original problem~\ref{eq:vp_sample}, this optimization contains no integer variables and is therefore computationally feasible for large numbers of nodes and training examples; it also naturally provides full bid curves. 
We note that electricity markets typically specify a maximum number of allowed segments in a bid curve.
We found that with reasonable parameters, the number of segments with non-zero volumes rarely exceeds 3 segments.
However, if required, the number of segments can be adjusted post-optimization; the simplest way is to discard the segments with the smallest absolute volume.

In the following sections, we show several special cases of the general framework~\ref{eq:general_so}, which apply some simplifications and approximations and encompass other convergence bidding strategies in the literature (Figure~\ref{fig:vp_v_p_models}).

\subsection{Neglecting the bid price: self-scheduling models (V-models)}

A common simplification of problem~\ref{eq:return_so} to obtain bid volumes is to use single-segment bids and set the supply and demand prices such that the bids are cleared almost surely.
Mathematically, this corresponds to $\overline{p} \to -\infty$ and $\underline{p} \to \infty$, which means that the day-ahead price variables can be marginalized out from the joint distribution.
For instance, in the case of the expected value of our specific definition of revenue, we can write:
\begin{equation}
    \begin{aligned}
        & \mathbb{E}_{p(\Delta, \lambda)}\left[(1_{\lambda \succeq -\infty} \odot \overline{w} + 1_{\lambda \preceq \infty} \odot   \underline{w})^{T} \Delta \right] = \\
        & \quad = \int \int (1_{\lambda \succeq -\infty} \odot \overline{w} + 1_{\lambda \preceq \infty} \odot   \underline{w})^{T} \Delta p(\Delta, \lambda) d\Delta d\lambda \\
        & \quad = \int \int (\overline{w} + \underline{w})^{T} \Delta p(\Delta, \lambda) d\Delta d\lambda \\
        & \quad = \int (\overline{w} + \underline{w})^{T} \Delta \left( \int p(\Delta, \lambda) d\lambda \right) d\Delta \\
        & \quad = \int (\overline{w} + \underline{w})^{T} \Delta p(\Delta) d\Delta \\
        & \quad = \mathbb{E}_{p(\Delta)}\left[ (\overline{w} + \underline{w})^{T} \Delta \right].
    \end{aligned}
\end{equation}
For other functionals (e.g. variance, value at risk, expected shortfall) $\lambda$ can be similarly marginalized out.
In this class of models, therefore, only $p(\Delta)$ needs to be modelled and only volumes are used as decision variables:
\begin{equation}
    \begin{aligned}
        \max_{\substack{{\overline{w} \in \mathbb{R}^{N}} \\ {\underline{w} \in \mathbb{R}^{N}}}} & \ \mathbb{F}_{p(\Delta)}\left[\left(\overline{w} + \underline{w} \right)^{T} \Delta \right] \\
        \mathrm{s.t.} & \ \mathbb{G}_{p(\Delta)}\left[\left(\overline{w} + \underline{w} \right)^{T} \Delta \right] \le 0, \\
        & \ (\overline{w}, \underline{w}) \in \mathcal{W}. \\
    \end{aligned}
    \label{eq:return_v}
\end{equation}

These approaches fall into the category of self-scheduling strategies identified in~\cite{Samani21}:
\begin{enumerate}
    \item The main focus is on the bid volumes and the bid prices are set to such values that guarantee the clearance of these bids (i.e. low and high values for supply and demand positions, respectively).
    \item $p(\Delta)$ is modelled and used to obtain optimal bid volumes and therefore there is not necessarily high correlation of the optimal position (sign of volume) among the subsequent submitted bids on the same node.
    \item Since the bid price is not modelled, the bids are single-step ($S = 1$).
\end{enumerate}

\subsubsection*{Traditional portfolio optimizations: Markowitz (mean-variance) and mean-expected shortfall}

Using a probabilistic model with a multivariate normal distribution $p(\Delta | x^{*}, X, Y, \theta) = \mathcal{N}(\Delta | \mu(x^{*}, X, Y, \theta), \Sigma(x^{*}, X, Y, \theta))$ with mean vector $\mu$ and covariance matrix $\Sigma$ to approximate $p(\Delta)$ and taking the variance of the revenue as risk measure leads to the well-known Markowitz portfolio optimization:
\begin{equation}
    \begin{aligned}
        \max_{\substack{{\overline{w}, \underline{w} \in \mathbb{R}^{N}}}} & \ \mathbb{E}_{\mathcal{N}(\Delta | \mu, \Sigma)}\left[(\overline{w} + \underline{w})^{T} \Delta\right] \\
        \mathrm{s.t.} & \ \mathbb{V}_{\mathcal{N}(\Delta | \mu, \Sigma)}\left[(\overline{w} + \underline{w})^{T} \Delta\right] \le \rho, \\
        & \ (\overline{w}, \underline{w}) \in \mathcal{W}, \\
    \end{aligned}
    \label{eq:v_anal_markowitz}
\end{equation}
where $\rho$ again denotes the upper bound on the risk.
It is easy to show that in the above problem, $\mathbb{E}_{\mathcal{N}(\Delta | \mu, \Sigma)}\left[(\overline{w} + \underline{w})^{T} \Delta\right] = (\overline{w} + \underline{w})^{T} \mu$ and $\mathbb{V}_{\mathcal{N}(\Delta | \mu, \Sigma)}\left[(\overline{w} + \underline{w})^{T} \Delta\right] = (\overline{w} + \underline{w})^{T} \Sigma (\overline{w} + \underline{w})$, so we can get the traditional Markowitz setup:
\begin{equation}
    \begin{aligned}
        \max_{\substack{{\overline{w}, \underline{w} \in \mathbb{R}^{N}}}} & \ (\overline{w} + \underline{w})^{T} \mu \\
        \mathrm{s.t.} & \ (\overline{w} + \underline{w})^{T} \Sigma (\overline{w} + \underline{w}) \le \rho, \\
        & \ (\overline{w}, \underline{w}) \in \mathcal{W}. \\
    \end{aligned}
    \label{eq:v_anal_markowitz_trad}
\end{equation}
The above approach was used for instance in \cite{Scott13}.
Similarly, one can use the expected shortfall instead of the variance as risk measure for a normal distribution:
\begin{equation}
    \begin{aligned}
        \max_{\substack{{\overline{w}, \underline{w} \in \mathbb{R}^{N}}}} & \ \mathbb{E}_{\mathcal{N}(\Delta | \mu, \Sigma)}\left[(\overline{w} + \underline{w})^{T} \Delta\right] \\
        \mathrm{s.t.} & \ \mathbb{ES}^{\alpha}_{\mathcal{N}(\Delta | \mu, \Sigma)}\left[(\overline{w} + \underline{w})^{T} \Delta\right] \le \rho, \\
        & \ (\overline{w}, \underline{w}) \in \mathcal{W}. \\
    \end{aligned}
    \label{eq:v_anal_es}
\end{equation}
For a normal distribution of delta prices, the expected shortfall of the revenue has a closed-from expression, which is proportional to the standard deviation of the revenue, i.e. $\sqrt{(\overline{w} + \underline{w})^{T} \Sigma (\overline{w} + \underline{w})}$.
Therefore, both the Markowitz and expected shortfall PO formulations are second-order cone programming (SOCP) problems.

\subsubsection*{Sample-based mean-expected shortfall portfolio optimization}

Alternatively, $p(\Delta)$ can be approximated directly by the observed delta prices.
In this case, a simple LP can be constructed:
\begin{equation}
    \begin{aligned}
        \max_{\substack{{\overline{w}, \underline{w} \in \mathbb{R}^{N}} \\ {z \in \mathbb{R}^{T}, \tau \in \mathbb{R}}}} & \ \frac{1}{T} \sum \limits_{t=1}^{T} r_t \\
        \mathrm{s.t.} & \ -\tau + \frac{1}{K} \sum \limits_{t=1}^{T} z_{t} \le \rho, \\
        & \ z_{t} \ge 0, & \forall t \\
        & \ z_{t} \ge \tau - r_{t}, & \forall t \\
        & \ r_{t} = (\overline{w} + \underline{w})^{T} \Delta_{t}, & \forall t \\
        & \ (\overline{w}, \underline{w}) \in \mathcal{W}, \\
    \end{aligned}
    \label{eq:v_sample_es}
\end{equation}
where $K = \lfloor \alpha T \rfloor$ as before. 
The above problem serves as the basis in~\cite{Wang19} and~\cite{Li22}, where instead of using continuous variables for the supply and demand volumes, binary variables representing 0 or 1 MWh of convergence bid are applied leading to a mixed-integer linear programming and mixed-integer quadratic programming problems, respectively.

\subsection{Neglecting the bid volume: opportunistic models (P-models)}

Another direction to simplify problem~\ref{eq:return_so} is to get rid of the volume variables and focus solely on the bid prices.
For instance, by setting the maximum number of segments to one and using constant nodal volumes, i.e. $\overline{w}^{n} = 1$ and $\underline{w}^{n} = -1$ for $n=1,\dots,N$ in problem \ref{eq:return_so} we have the following optimization problem:
\begin{equation}
    \begin{aligned}
        \max_{\substack{{\overline{p} \in \mathbb{R}^{S \times N}} \\ {\underline{p} \in \mathbb{R}^{S \times N}}}} & \ \mathbb{F}_{p(\Delta, \lambda)}\left[\left(1_{\lambda \succeq \overline{p}} - 1_{\lambda \preceq \underline{p}} \right)^{T} \Delta\right] \\
        \mathrm{s.t.} & \ \mathbb{G}_{p(\Delta, \lambda)}\left[\left(1_{\lambda \succeq \overline{p}} - 1_{\lambda \preceq \underline{p}} \right)^{T} \Delta \right] \le 0, \\
        & \ (\overline{p}, \underline{p}) \in \mathcal{P}. \\
    \end{aligned}
    \label{eq:return_p}
\end{equation}

It is easy to show that the expected value of the revenue expression in the above optimization problem becomes the sum of the expected values of the nodal revenues of the supply and demand positions:
\begin{equation}
    \begin{aligned}
        & \mathbb{E}_{p(\Delta, \lambda)}\left[\left(1_{\lambda \succeq \overline{p}} - 1_{\lambda \preceq \underline{p}} \right)^{T} \Delta \right] = \\
        & \quad = \mathbb{E}_{p(\Delta, \lambda)}\left[ \sum_{n=1}^{N} \left(1_{\lambda^{n} \ge {\overline{p}}^{n}} - 1_{\lambda^{n} \le {\underline{p}}^{n}} \right) \Delta^{n} \right] \\
        & \quad = \sum_{n=1}^{N} \mathbb{E}_{p(\Delta, \lambda)}\left[\left(1_{\lambda^{n} \ge {\overline{p}}^{n}} - 1_{\lambda^{n} \le {\underline{p}}^{n}} \right) \Delta^{n} \right] \\
        & \quad = \sum_{n=1}^{N} \mathbb{E}_{p(\Delta^{n}, \lambda^{n})}\left[\left(1_{\lambda^{n} \ge {\overline{p}}^{n}} - 1_{\lambda^{n} \le {\underline{p}}^{n}} \right) \Delta^{n} \right] \\
        & \quad = \sum_{n=1}^{N} \mathbb{E}_{p(\Delta^{n}, \lambda^{n})}\left[1_{\lambda^{n} \ge {\overline{p}}^{n}} \Delta^{n} \right] - \mathbb{E}_{p(\Delta^{n}, \lambda^{n})}\left[1_{\lambda^{n} \le {\underline{p}}^{n}} \Delta^{n} \right],
    \end{aligned}
\end{equation}
where $p(\Delta^{n}, \lambda^{n})$ for $n=1,\dots,N$ are the nodal marginal distributions.
However, for an arbitrary functional, the revenue expression does not necessarily simplify to separate expressions of nodal contributions and positions due to the correlation between the day-ahead and delta price components, which makes it challenging to derive a closed-form optimization problem.

\subsubsection*{Sample-based mean-expected shortfall price model}
Assuming zero correlation between different nodal variables (i.e. $p(\Delta, \lambda) = \prod \limits_{n=1}^{N} p(\Delta^{n}, \lambda^{n})$) and approximating the $p(\Delta^{n}, \lambda^{n})$ distributions by the observed prices lead to a variant of the dynamic node labelling approach~\cite{Saez-Gallego18, Samani21}, where for each node and each position, a separate optimization problem is solved.
For instance, for node $n$ and supply position we can write:
\begin{equation}
    \begin{aligned}
        \max_{\substack{{\overline{p}^{n} \in \mathbb{R}} \\ {z \in \mathbb{R}^{T}, \tau \in \mathbb{R}}}} & \ \frac{1}{T}\sum \limits_{t=1}^{T} r_{t} \\
        \mathrm{s.t.} & \ -\tau + \frac{1}{K} \sum \limits_{t=1}^{T} z_{t} \le \rho, \\
        & \ z_{t} \ge 0, & \forall t \\
        & \ z_{t} \ge \tau - r_{t}, & \forall t \\
        & \ r_{t} = \begin{cases} \Delta_{t}^{n} & \text{if } \lambda_{t}^{n} \ge \overline{p}^{n} \\ 0 & \text{else} \end{cases}, & \forall t \\
        & \ \overline{p}_{\mathrm{min}}^{n} \le \overline{p}^{n} \le \overline{p}_{\mathrm{max}}^{n}. \\
    \end{aligned}
    \label{eq:dnl_es_supply}
\end{equation}
Note that in the above problem - similarly to problem \ref{eq:vp_sample} - binary variables are required leading to a MILP optimization problem.

The above model illustrates that P-models tend to fall into the category of opportunistic strategy identified by~\cite{Samani21}:
\begin{enumerate}
    \item The main focus is on the bid prices and bid volumes are set by some simple heuristics.
    \item By setting the bid price, these models are able to conditionally clear in favourable situations such as price spikes, and therefore the optimal bid prices might tend to be rejection-prone (i.e. high and low values for supply and demand positions, respectively) leading to high correlation of the positions among the subsequent submitted bids.
    \item Since multiple optimal bid prices are also possible, bid curves with multiple segments can be also obtained. 
\end{enumerate}

\subsubsection*{Sample-based mean-expected shortfall bid curve price model}

Similarly to the sample-based volume-price model, the cleared-rejected delta prices can be precomputed.
Also, by introducing a weight vector $\overline{w}^{n} \in \mathbb{R}^{\overline{P}}$ for the allowed supply bid prices (again, taken from the training day-ahead prices) for node $n$ and distribute the unit volume among these as bid segments, we can construct the following linear programming problem:
\begin{equation}
    \begin{aligned}
        \max_{\substack{{\overline{w}^{n} \in \mathbb{R}^{\overline{P}}} \\ {z \in \mathbb{R}^{T}, \tau \in \mathbb{R}}}} & \ \frac{1}{T} \sum \limits_{t=1}^{T} r_{t} \\
        \mathrm{s.t.} & \ -\tau + \frac{1}{K} \sum \limits_{t=1}^{T} z_{t} \le \rho, \\
        & \ z_{t} \ge 0, & \forall t \\
        & \ z_{t} \ge \tau - r_{t}, & \forall t \\
        & \ r_{t} = \overline{\Delta}^{n} \overline{w}^{n}, & \forall t \\
        & \ 1^{T} \overline{w}^{n} \le 1, \\
        & \ \overline{w} \succeq 0, \\
        \mathrm{with} &\ \overline{\Delta}^{n} = (\Delta^{n} \otimes 1_{\overline{P}}) \odot (\lambda^{n} \otimes_{\ge} \overline{p}^{n}), \\
        & \ \overline{p}^{n} \in \overline{\mathcal{P}}^{n}.
    \end{aligned}
    \label{eq:fast_bc_dnl_es_sample}
\end{equation}

The optimization problem for the demand position can be similarly constructed.

\subsection{Technical details of numerical experiments}

In the previous sections, we introduced the general form of VP-models (problem~\ref{eq:general_so}) that use both the bid volumes and prices as decision variables and we also discussed its simplifications: V-models (problem~\ref{eq:return_v}), where only bid volumes are considered as decision variables and P-models (problem~\ref{eq:return_p}), where only bid prices are used as decision variables.
We also specified the corresponding LP-based formulations using training samples directly.
Among these, according to our best knowledge, sample-VP (\ref{eq:fvp_sample}) is the first concrete model in the literature that tries to produce optimal bid volumes and prices simultaneously based on some training data, while
the sample-V (\ref{eq:v_sample_es}) model is a sample-based version of the widely used traditional portfolio optimization approaches~\cite{Scott13, Samani21} and the sample-P (\ref{eq:fast_bc_dnl_es_sample}) model is a fast version of the dynamic node labelling approach~\cite{Saez-Gallego18, Samani21} that uses expected shortfall as risk measure instead of the revenue-loss ratio.

In order to demonstrate the advantage of obtaining bid volumes and prices simultaneously over the volume- and price-based models, we performed numerical experiments using the sample-based models and compared their financial performance using real-world price data between 2017-2021 from the California ISO~\cite{OASIS}.

We note that in general it is rather challenging to compare the performance of state-of-the-art models as they use different optimization formulations, risk measures, corresponding upper bound or risk aversion parameters, training data size and predicted distribution of the joint day-ahead and delta prices.
Our common framework, however, makes it possible to perform a fair comparison between the above-mentioned models using the same formulation (i.e. maximizing the expected value of revenues with an upper bound on the risk), risk measure (i.e. expected shortfall of revenues) and training data.
We also emphasize that focusing solely on the achieved cumulative revenue within some test period might be misleading as it does not take potentially large losses into account.
Therefore, when comparing the performance of the different models, we investigate the trade-off between the expected value and expected shortfall of the produced revenues.
This comparison strategy is also aligned with the formulation that these optimization models use.

The effective number of nodes was reduced by using an event synchronization-based correlation on the biddable nodes of CAISO~\cite{Cui18} with a threshold of 0.98.
Depending on the target day, this resulted in approximately 500 - 750 nodes.

For all models, an hourly optimization with a 1-year rolling-window was applied, i.e. for each target hour, an independent optimization was performed using 1-year of training data of the same hour of day delta and day-ahead prices (365 samples per node), which resulted in simulations over a 4-year long period (2018-2021). 

The maximum expected revenue with upper bound on the expected shortfall formulation was used for all optimization models (i.e. problems \ref{eq:fvp_sample}, \ref{eq:v_sample_es} and \ref{eq:fast_bc_dnl_es_sample} for sample-VP, sample-V and sample-P models, respectively). 
P-models do not include volumes as decision variables (i.e. they apply unit volumes), therefore, in order to have a consistent risk measure in all models, we introduced a volume-normalized expected shortfall, $\tilde{\rho}$, measured in \$/MWh and the actual expected shortfall in the corresponding problems can be obtained by multiplying it by the total volume used in the optimization model ($\rho = W \tilde{\rho}$).
This guarantees a fair comparison between models using joint nodal-based optimization (sample-VP and sample-V models) and single nodal-based optimization (sample-P models).
For computing the expected shortfall, we used a 5\% quantile-level ($\alpha = 0.05$).

For each hour, a maximum of 1 GWh volume was distributed among the nodes.
For the sample-P models, the most promising 100 supply and 100 demand nodes were selected (i.e. producing the highest expected revenue values by the corresponding model) and a maximum of 5 MWh volume was allocated for each (note that the sample-P model provides a normalized volume, which is then multiplied by 5).
Note that there could be overlap among the nodes appearing in the most promising supply and demand sets. 
In order to have a fair comparison, the same 200 positions were used in the sample-VP and sample-V models as well with a maximum nodal volume of 50 MWh.
For sample-V models, the corresponding supply and demand prices were chosen as the market caps, to ensure a high clearance rate.
For the sample-VP and sample-P models, there was not any particular constraint on the bid prices. 
The minimum submittable volume allowed in CAISO is 1 MWh, so we simply filtered out bid segments with less volume than this.
Also, the maximum number of bid segments is 10 in CAISO, so in the very rare cases when the number of derived bid segments exceeded the maximum allowed number, we simply removed excess segments with the lowest volume.

We emphasize that the models we deal in these experiments do not include any consideration of the impacts of submitted bids on the dayahead market.
We also note that fees were not included in these experiments and statistics were computed on gross revenue only.
Simulations were performed within the JuMP package \cite{Dunning17} written in Julia \cite{Bezanson12} using the Clp  \cite{Forrest22clp} and Cbc \cite{Forrest22cbc} solvers.

\section{Numerical Results}
\label{sec:results}

The volume-normalized expected shortfall ($\tilde{\rho}$) makes it possible to perform a fair comparison among the 3 main types of models we consider.
Following \cite{Samani21}, we used 1-year long training data in these models, so the hyperparameter we changed and investigated first was the normalized risk.
We performed simulations for the sample-based P-, V- and VP-models with a varying upper limit of the normalized-expected shortfall, using 0.1, 1 and 10 \$/MWh values.

While in the sample-P model, the 1 GWh hourly volume is evenly distributed among the most promising 100 supply and 100 demand positions, the sample-V and sample-VP models are proper portfolio optimization approaches that place more volume on the more promising positions (up to the maximum nodal volume of 50 MWh) over the less beneficial ones.
Therefore, we expect these models to have a better performance in general.
In order to improve the performance of the sample-P model, we also investigated another case, when the 1 GWh was distributed only among the most beneficial 10 supply and 10 demand positions with the maximal 50 MWh nodal volume (model sample-P-max).

\begin{figure}[!ht]
    \centerline{\includegraphics[width=\columnwidth]{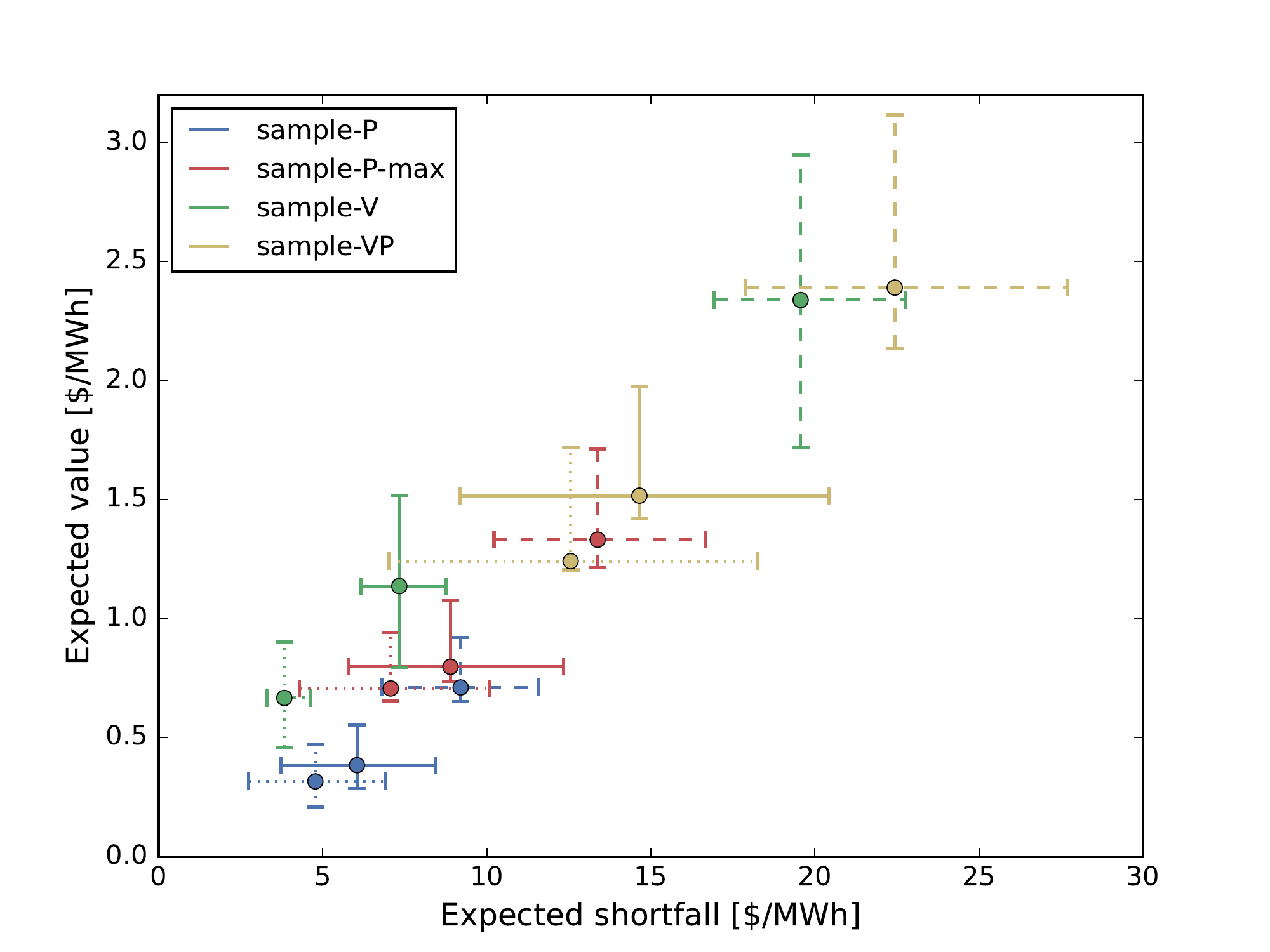}}
    \caption{Two-sided confidence intervals with 95\% confidence level of expected value and expected shortfall of the hourly volume-normalized revenues of the investigated models using 0.1 (dotted), 1 (solid) and 10 (dashed) values for the upper limit of the volume-normalized expected shortfall with a 1-year long training data size.}
    \label{fig:ci}
\end{figure}

For all models, we computed the expected revenue (higher is better), $\int \tilde{r} p(\tilde{r}) d\tilde{r}$, as well as the expected shortfall (lower is better), $-\frac{1}{\alpha} \int_{\tilde{r} \le q_{\alpha}(\tilde{r})}\tilde{r} p(\tilde{r}) d\tilde{r}$ and windfall (higher is better), $\frac{1}{\alpha} \int_{\tilde{r} \ge q_{1-\alpha}(\tilde{r})} \tilde{r} p(\tilde{r}) d\tilde{r}$ values of the hourly volume-normalized revenues ($\tilde{r}$) for the simulation period (where $q_{\alpha}$ is the quantile function with $\alpha$ quantile-level set to $0.05$ as in the optimization models).
The results are collected in Table~\ref{tab:return_statistics}.
In general, the sample-P model has the lowest hourly expected revenue.
Using the maximum nodal volumes for the best 20 positions (sample-P-max) significantly improves the performance and the corresponding sample-P-max model provides almost twice the expected revenues.
The V- and VP-models show significantly higher expected values and the sample-VP model has the highest ones for each normalized risk parameter value.
The corresponding expected shortfall and expected windfall values show a similar trend indicating that models with higher expected revenue have a higher risk as well.
This observation is further supported by Figure~\ref{fig:ci}, which shows the corresponding two-sided confidence intervals of the expected value and expected shortfall of the hourly normalized revenues with 95\% confidence level.
The confidence intervals were estimated by a simple subsampling approach with a 2-year long moving block size \cite{Politis99}.

The significant difference in the performance among the approaches is likely due to the difference of how these models treat the mean-risk trade-off: while the sample-P models consider only the joint delta and day-ahead price distribution of each single node separately ($p(\Delta^{n}, \lambda^{n})$), the sample-V model uses the joint delta price distribution ($p(\Delta)$) and the sample-VP approach uses the joint delta and day-ahead price distribution of all positions ($p(\Delta, \lambda)$) to compute the hourly revenue values in their optimization model.
Considering the nodal price distributions together makes it possible to significantly increase expected revenue by distributing the available total volume among the more beneficial nodes and positions.
However, the higher dimensionality of the distributions of the sample-V and sample-VP models has two important consequences.
First, these models are more exposed to non-stationarity.
And second, they can obtain less accurate revenue statistics with the same sample size.
These two effects result in higher expected shortfall. 

\begin{table}[!h]
    \caption{Expected value, shortfall and windfall statistics of the hourly normalized revenues in the 4-year long simulation period for the different models.}
    \label{tab:return_statistics}
    \centering
    \begin{tabular}{lrrrr}
        \toprule
        \multirow{2}{*}{Model} & $\tilde{\rho}$ & \multicolumn{3}{c}{Hourly revenue statistics [\$/MWh]} \\
        \cmidrule(r){3-5}
        & [\$/MWh] & Exp. value & Exp. shortfall & Exp. windfall \\
        \midrule
        sample-P & 0.1 & 0.317 & 4.775 & 9.362 \\
         & 1 & 0.385 & 6.041 & 10.830 \\
         & 10 & 0.712 & 9.202 & 13.832 \\
        \midrule
        sample-P-max & 0.1 & 0.707 & 7.071 & 19.272 \\
         & 1 & 0.798 & 8.892 & 20.908 \\
         & 10 & 1.332 & 13.379 & 27.354 \\
        \midrule
        sample-V & 0.1 & 0.667 & 3.829 & 11.523 \\
         & 1 & 1.137 & 7.332 & 19.563 \\
         & 10 & 2.339 & 19.564 & 45.255 \\
        \midrule
        sample-VP & 0.1 & 1.242 & 12.554 & 29.207 \\
         & 1 & 1.517 & 14.650 & 31.659 \\
         & 10 & 2.392 & 22.433 & 49.762 \\
        \bottomrule
    \end{tabular}
\end{table}

\begin{figure}[!ht]
    \centerline{\includegraphics[width=\columnwidth]{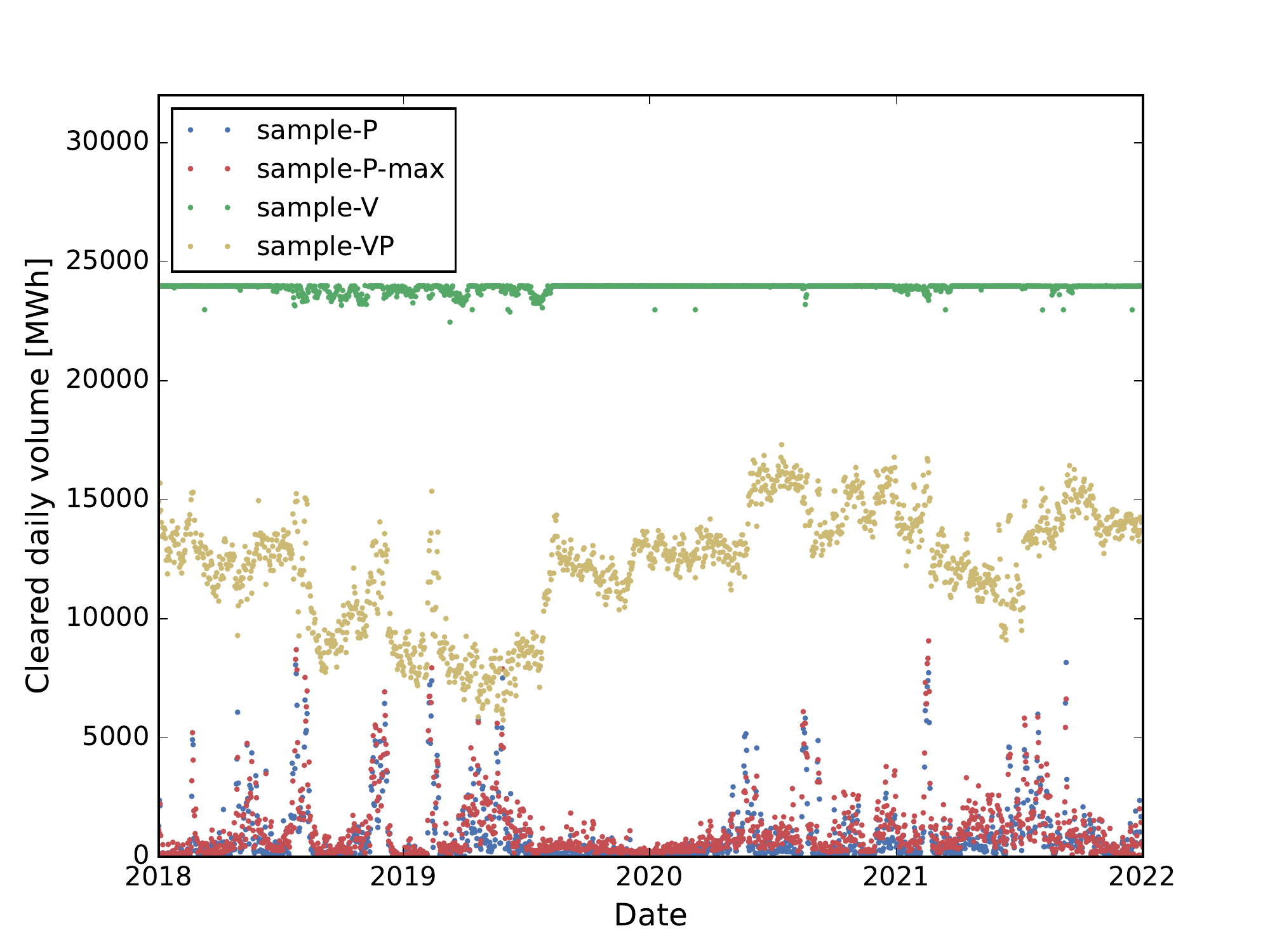}}
    \caption{Daily cleared volumes for the models with $\tilde{\rho} = 1$ \$/MWh normalized risk parameter.}
    \label{fig:volume}
\end{figure}

Although the attempted hourly volumes are set to approximately 1 GWh for each model, the cleared volumes can differ significantly.
For instance, Figure~\ref{fig:volume} shows the daily aggregated cleared volumes for the models with $\tilde{\rho} = 1$ \$/MWh normalized risk.
Not surprisingly, the sample-V model has a clearance rate close to 100\% (24 GWh per day), given the corresponding supply and demand bid prices are rather low and high, respectively.
For sample-VP, the cleared daily volume is significantly smaller and most of the time it does not exceed 15 GWh.
Finally, as expected, both sample-P models are rather rejection-prone approaches, whose cleared daily volume rarely exceeds 5 GWh.
This is in accordance with the fact that they are opportunistic models with a high rejection rate.

In Table~\ref{tab:volume_statistics}, we collected the hourly volume statistics for the different models including the attempted and cleared total volumes as well as their supply and demand rates.
It is clear from the table that when increasing the risk parameter, the cleared volume increases in all models.
For the sample-V model, it already has a high clearance rate even with $\tilde{\rho} = 0.1$ \$/MWh (note that although in its model the corresponding supply and bid prices are assumed to be $-\infty$ and $\infty$ respectively, in the actual bids these numbers were finite values).
The sample-P model has the lowest hourly average cleared volumes with only 17, 34 and 303 MWh for 0.1, 1 and 10 \$/MWh values of the risk parameter, respectively.
Although its financial performance is significantly better, the sample-P-max model has just slightly larger average cleared volumes.
As was already indicated in Figure~\ref{fig:volume}, the most promising sample-VP approach obtains its highest revenue with significantly lower cleared volumes compared to the sample-V model for all risk parameters.
For example, for $\tilde{\rho} = 10$ \$/MWh, where their expected revenue is comparable (cf. Table~\ref{tab:return_statistics}), the sample-VP model achieves it by an average cleared volume of 74.2\% of that of the sample-V model, which further favors the volume-price model if fees are also considered.

Looking at the distribution of volumes between supply and demand positions, it seems that for the sample-V model, it is close to 50-50\% for both the attempted and cleared volumes.
For sample-P models, where the attempted supply and demand volumes were set to equal, in the cleared volumes there is a larger portion of supply positions.
Finally, although the sample-VP model preferred supply positions, the cleared volumes were more balanced between the supply and demand positions. 

\begin{table*}[!t]
    \caption{Hourly volume statistics of the investigated models.}
    \label{tab:volume_statistics}
    \centering
    \begin{tabular}{lrrrrr}
        \toprule
        \multirow{2}{*}{Model} & $\tilde{\rho}$ & \multicolumn{2}{c}{Attempted volume} & \multicolumn{2}{c}{Cleared volume}\\
        \cmidrule(r){3-6}
        & [\$/MWh] & Mean total [\$/MWh] & Supply / demand [\%] & Mean total [\$/MWh] & Supply / demand [\%] \\
        \midrule
        sample-P & 0.1 & 938 & 51.3 / 48.7 & 17 & 69.6 / 30.4 \\
        & 1 & 911 & 51.6 / 48.4 & 34 & 64.5 / 35.5 \\
        & 10 & 898 & 50.7 / 49.3 & 303 & 59.0 / 41.0 \\
        \midrule
        sample-P-max & 0.1 & 996 & 50.0 / 50.0 & 22 & 63.4 / 36.6 \\
        & 1 & 992 & 50.1 / 49.9 & 46 & 62.3 / 37.7 \\
        & 10 & 988 & 50.3 / 49.7 & 335 & 57.9 / 42.1 \\
        \midrule
        sample-V & 0.1 & 986 & 50.2 / 49.8 & 941 & 50.2 / 49.8 \\
         & 1 & 998 & 51.6 / 48.4 & 996 & 51.6 / 48.4 \\
         & 10 & 1000 & 65.2 / 34.8 & 1000 & 65.2 / 34.8 \\
        \midrule
        sample-VP & 0.1 & 996 & 72.8 / 27.2 & 324 & 56.0 / 44.0 \\
         & 1 & 999 & 66.8 / 33.2 & 512 & 56.1 / 43.9 \\
         & 10 & 1000 & 68.6 / 31.4 & 742 & 66.8 / 33.2 \\
        \bottomrule
    \end{tabular}
\end{table*}

A comparison of bid statistics from the models is shown in Table~\ref{tab:bid_statistics}.
First, we investigated what fractions of the submitted bids are single position (i.e. for a given hour and node, the bid curve is either a supply or a demand position), or double position (i.e. for a given hour and node, both supply and demand bid curves were submitted).
We found that the sample-P model has the highest fraction of submitting double position bid curves, while the sample-VP and sample-P-max models have significantly smaller fractions (for sample-V no double position bids were allowed).

For models that produce multistep bid curves, we found that the highest number of segments is provided by the sample-P-max model, which reached the maximum number of 10 segments for larger normalized risk parameters.
The maximum number of segments was much smaller for the sample-P model: the corresponding segments were filtered out due to their small volume.
Interestingly, the maximum number of segments decreased for the sample-VP model with the increase of the risk parameter.

We also computed the fraction of bid curves with different number of segments.
For sample-P models, the fraction of multistep bids increases with increasing risk parameter.
For the sample-VP model, however, single-step bids dominate for higher risk parameters.

\begin{table*}[!t]
    \caption{Statistics of the derived bids of the investigated models.}
    \label{tab:bid_statistics}
    \centering
    \begin{tabular}{lrrrrrrr}
        \toprule
        \multirow{2}{*}{Model} & $\tilde{\rho}$ & \multicolumn{2}{c}{Fraction of positions [\%]} & \multirow{2}{*}{Maximum \# of segments} & \multicolumn{3}{c}{Fraction of segments [\%]} \\
        \cmidrule{3-4} \cmidrule{6-8}
        & [\$/MWh] & Single & Double & & Single-step & Double-step & $>$Double-step \\
        \midrule
        sample-P & 0.1 & 85.9 & 14.1 & 3 & 85.2 & 14.8 & 0.001 \\
        & 1 & 85.7 & 14.3 & 4 & 74.2 & 25.8 & 0.02 \\
        & 10 & 87.3 & 12.7 & 4 & 58.8 & 34.8 & 6.4 \\
        \midrule
        sample-P-max & 0.1 & 97.9 & 2.1 & 4 & 53.0 & 47.0 & 0.02 \\
        & 1 & 97.8 & 2.2 & 10 & 31.0 & 60.8 & 8.2 \\
        & 10 & 96.4 & 3.6 & 10 & 31.5 & 32.8 & 35.7 \\
        \midrule
        sample-V & 0.1 & 100.0 & 0.0 & 1 & 100.0 & 0.0 & 0.0 \\
         & 1 & 100.0 & 0.0 & 1 & 100.0 & 0.0 & 0.0 \\
         & 10 & 100.0 & 0.0 & 1 & 100.0 & 0.0 & 0.0 \\
        \midrule
        sample-VP & 0.1 & 96.2 & 3.8 & 9 & 79.6 & 15.4 & 5.0 \\
         & 1 & 96.1 & 3.9 & 8 & 79.5 & 16.3 & 4.2 \\
         & 10 & 97.9 & 2.1 & 5 & 93.8 & 5.7 & 0.5 \\
        \bottomrule
    \end{tabular}
\end{table*}

So far, we investigated the behaviour of the models by varying the normalized risk parameter with a fixed 1-year long training data size.
We also examined the effect of using shorter training windows of half-year (180 days) and quarter-year (90 days) periods with a fixed risk parameter $\tilde{\rho} = 1$ \$/MWh.
In principle, reducing the training data size in these sample-based models has a mixed effect.
As discussed earlier, on one hand, shorter training windows reduce the exposure to the non-stationary behaviour of $p(\Delta, \lambda)$.
On the other hand, using a smaller number of samples increases the statistical error of the computed expected value and risk of revenues.
Given the effective number of samples used to compute the expected shortfall is just a fraction of the total number of samples (e.g. 5\% for $\alpha = 0.05$), this quantity is more sensitive to this effect.

\begin{figure}[!ht]
    \centerline{\includegraphics[width=\columnwidth]{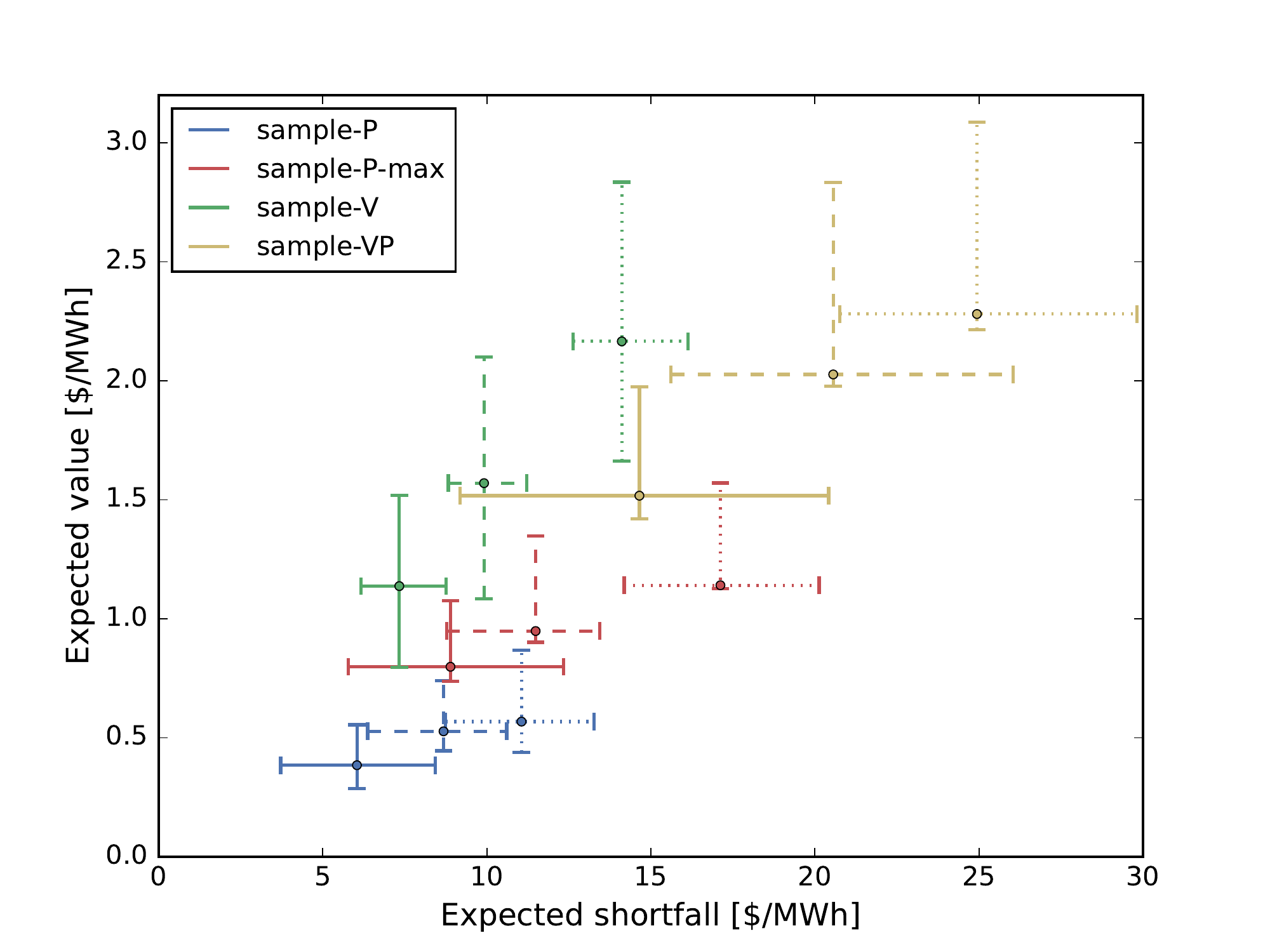}}
    \caption{Two-sided confidence intervals with 95\% confidence level of expected value and expected shortfall of the hourly volume-normalized revenues of the investigated models using quarter-year (dotted), half-year (dashed) and 1-year (solid) long training data sizes with $\tilde{\rho} = 1$ \$/MWh normalized risk parameter.}
    \label{fig:ci_lb}
\end{figure}

Figure~\ref{fig:ci_lb} shows the confidence intervals of the expected value and expected shortfall of revenues of the models using different training data sizes (confidence intervals were computed similarly as before).
It is clear from the figures that reducing the training data size from 1-year increases the total (expected value of) revenue.
The extent of increase is more significant for the sample-V and sample-VP models.
However, at the same time, the expected shortfall also increased for all approaches and models including volumes as decision variables are again more affected.
These observations indicate that there indeed exists a trade-off between the effects of non-stationarity and statistical accuracy.

\begin{figure}[!ht]
    \centerline{\includegraphics[width=\columnwidth]{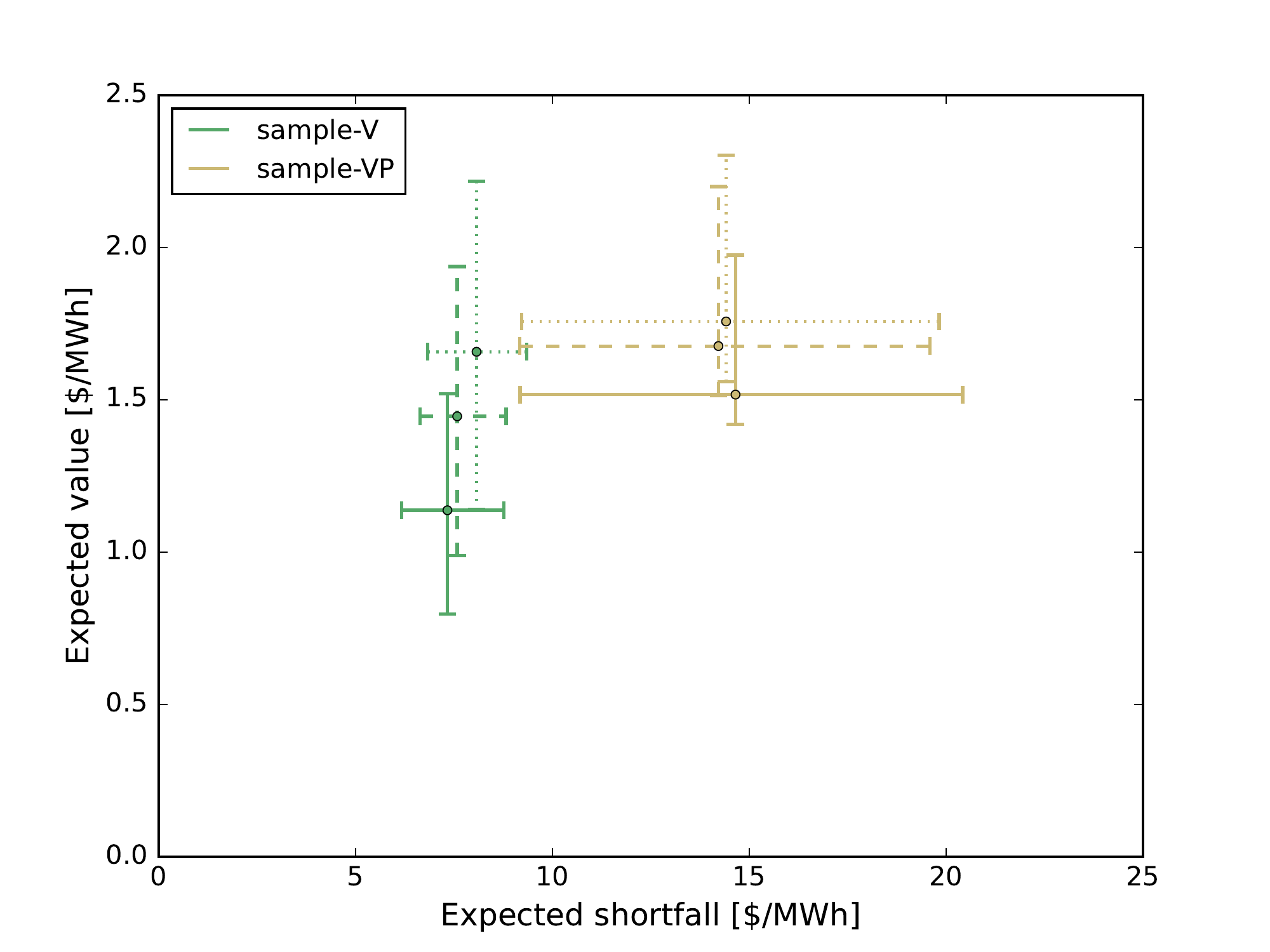}}
    \caption{Two-sided confidence intervals with 95\% confidence level of expected value and expected shortfall of the hourly volume-normalized revenues of V- and VP-models using the most promising 200 (solid) and the 400 (dashed), as well as all positions (dotted) with $\tilde{\rho} = 1$ \$/MWh normalized risk parameter.}
    \label{fig:ci_node}
\end{figure}

We also investigated the effect on the sample-V and sample-VP models when using a larger number of potential positions.
Figure~\ref{fig:ci_node} shows the comparison of expected values and expected shortfalls with confidence intervals using the most promising 200, 400 and all positions.
The positions were selected in the same way as before, i.e. based on the objective value of the sample-P model.
Using a larger number of positions significant increases the expected value of the revenue for the sample-V model.
Interestingly, for the sample-VP model, this improvement is much smaller, indicating that using a sample-P model for pre-selecting the most promising positions is a promising approach.
For the sample-V model, the improvement of the expected value is accompanied by a slight increase of the expected shortfall, while for the sample-VP model, the risk is basically unchanged.

Finally, we also compared the average computational times to solve the different optimization approaches for a target day (24 target hours) using a limited 200 positions, 400 positions and the full set of nodes, $2 \times 750$ nodes = 1500 positions, which was one of the largest problems in our experiments.
Optimizations were performed on a MacBook Pro (14-inch, 2021) machine with Apple M1 Max chip.
The results are collected in Table~\ref{tab:solve_times} for different risk parameters.
We begin with the comparison of P-models that use a MILP (problem~\ref{eq:dnl_es_supply}) and our LP (problem~\ref{eq:fast_bc_dnl_es_sample}) formulations.
It is clear that the mixed-integer linear programming based approach requires about an order of magnitude longer solve times than the corresponding linear programming ones.
This also demonstrates that using a MILP-based model for a VP-approach would be computationally extremely demanding for large $N$ and/or $T$.
Regarding the LP-based formulations, although sample-P models have the smallest optimization problems, a     large number of them needs to be solved (for each node and for each position).
Sample-V models show negligible solve times even for the largest problems.
Sample-VP models have the computationally most expensive optimization models but their solve times are still feasible even for the largest problems.

\begin{table}[!h]
    \addtolength{\tabcolsep}{-1.5pt}
    \caption{Average solve times of models with different number of positions.}
    \label{tab:solve_times}
    \centering
    \begin{tabular}{lrrrr}
        \toprule
        \multirow{2}{*}{Model} & $\tilde{\rho}$ & \multicolumn{3}{c}{Solve time [s]} \\
        \cmidrule(r){3-5}
        & [\$/MWh] & 200 pos. & 400 pos. & 1500 pos. \\
        \midrule
        sample-P & 0.1 & $129.0 \pm 5.1$ & $253.9 \pm 9.5$ & $930.1 \pm 40.7$ \\
        (MILP) & 1 & $139.3 \pm 5.0$ & $284.1 \pm 10.7$ & $1076.9 \pm 49.8$ \\
        & 10 & $145.9 \pm 16.8$ & $289.2 \pm 37.5$ & $1075.2 \pm 117.3$ \\
        \midrule
        sample-P & 0.1 & $15.0 \pm 0.1$ & $33.5 \pm 0.7$ & $119.2 \pm 0.6$\\
         & 1 & $16.2 \pm 0.3$ & $33.8 \pm 1.2$ & $113.1 \pm 1.9$ \\
         & 10 & $15.3 \pm 0.1$ & $31.1 \pm 0.4$ & $119.1 \pm 2.7$\\
        \midrule
        sample-V & 0.1 & $0.3 \pm 0.1$ & $0.7 \pm 0.1$ & $2.9 \pm 0.1$ \\
         & 1 & $0.3 \pm 0.1$ & $0.7 \pm 0.1$ & $2.9 \pm 0.1$ \\
         & 10 & $0.3 \pm 0.1$ & $0.6 \pm 0.1$ & $2.9 \pm 0.1$ \\
        \midrule
        sample-VP & 0.1 & $38.5 \pm 3.1$ & $103.0 \pm 9.3$ & $706.2 \pm 67.0$ \\
         & 1 & $37.6 \pm 2.5$ & $98.4 \pm 9.5$ & $637.6 \pm 43.4$ \\
         & 10 & $32.0 \pm 1.8$ & $79.1 \pm 5.3$ & $658.7 \pm 51.1$ \\
        \bottomrule
    \end{tabular}
\end{table}

\section{Conclusions}

Convergence bidding is an important element of efficient operation and scheduling of electric power markets in the United States.
The primary role of convergence bids submitted to the day-ahead market is to reduce the price gap between the day-ahead and real-time markets in two-settlement electricity market systems.
In this paper, we have introduced a general stochastic optimization framework for convergence bidding that provides both optimal supply and demand bid curves.
The general model includes both the bid volumes and prices as decision variables, and therefore we call this approach a volume-price model (VP-model).

The key object in our framework is the joint distribution $p(\Delta, \lambda)$ of the delta ($\Delta$) and day-ahead ($\lambda$) prices of all biddable nodes.
We discussed the advantages and disadvantages of approaches that approximate this distribution by a parametric model or by samples (either using empirical samples, or sampling a parametric model).
With appropriate assumptions, the framework simplifies to many approaches in the literature as special cases.
In the case of V-models, only bid volumes are considered as decision variables, while bid prices are usually set to such values that the bid is cleared with high probability (self-scheduling strategies).
For P-models, however, only bid prices are subject to optimize, while the corresponding volumes are constant (opportunistic strategies).

We showed that when using the sample-based approach with a formulation that maximizes the expected value of revenue with an upper bound on the expected shortfall of revenue, a linear program exists that can be efficiently solved to obtain optimal bid curves even for large markets.
This has several advantages over previous approaches: it does not require integer variables, it obtains co-optimized bid volumes and prices, and it naturally produces bid curves for both supply and demand positions.

Experiments with CAISO data indicate the improvement gained with this approach over approaches which optimize solely the bid volumes or the bid prices.
Performing a portfolio optimization that derives volumes based on either the joint delta price distribution (sample-V model) or joint delta and day-ahead price distribution of all considered positions (sample-VP model) improved performance compared to models using fixed volumes on the most promising positions (sample-P models).
We also found that simultaneously obtaining bid prices and volumes in sample-VP models outperformed the corresponding sample-V models with significantly less cleared volume.

The increased expected revenue was accompanied by an increased risk (expected shortfall).
Further investigations of varying training data size revealed that there is a trade-off between non-stationary effects and accuracy of computed statistics of revenues in these sample-based models.
Since in this paper our main purpose was to carry out a consistent comparison between different approaches, we did not perform extensive hyperparameter optimization.
However, we think that optimal hyperparameters (e.g. training data size, upper bound of risk etc.) could significantly alter the expected value--expected shortfall trade-off of the normalized revenue.

The general framework allows for a large degree of freedom in the representation of~$p(\Delta, \lambda)$.
It would be particularly interesting to investigate formulations where the joint distribution is represented by a parametric or semi-parametric model.
Adding an inductive bias in this manner can help to compensate for limited data availability, and thus help tackle the non-stationarity in~$p(\Delta, \lambda)$, which was demonstrated in our experiments.
Including exogenous data features to specify a more accurate~$p(\Delta, \lambda)$ is another natural extension.

Finally, we note that both the general framework and specific models described might be extended by several ways.
For example, transaction fees and uplift costs (cf.~\cite{Li22}) can be straightforwardly modelled by additional linear terms.
Also, distributionally-robust optimization approaches could be used to model the ambiguity set of distributions that are within a given probability metric from the empirical probability distribution (cf.~\cite{Du21}).
We emphasize that the models discussed here are fully price-based and use solely historical data.
Although models that treat bid prices as decision variables (P- and VP-models) can still benefit from unprecedented or infrequent price spikes, in general, they do not necessarily capture these rare events.
There are two potential directions, however, that can improve our framework for such conditions.
One way is to include exogenous variables in the models that can indicate price spikes and modify the training price distributions accordingly.
Another, more sophisticated direction is to model the day-ahead clearance problem as well in a bilevel optimization problem, along the lines suggested in section~\ref{sec:general_so}.
The additional advantage of such model is that the price impact of the submitted bids is also explicitly taken into account.
However, this approach has some computational challenges:
the resulting optimization problem would no longer be an LP due to the addition of the unit commitment problem, although it might be possible to integrate the sample-VP optimization with only a small overhead.
Also, modelling the behavior of other virtuals requires extensive sampling and solving a large number of optimization problems.

\section*{Acknowledgements}

The authors would like to thank Zoubin Ghahramani and Cozmin Ududec for fruitful discussions and their constructive and insightful comments, and Glenn Moynihan for invaluable contributions to software implementations.

\ifCLASSOPTIONcaptionsoff
  \newpage
\fi

\bibliographystyle{IEEEtran}

\bibliography{refs.bib}

\end{document}